\documentclass[12pt]{article}
\usepackage{amssymb,amsmath, amsthm}
\usepackage{graphicx}

\newtheorem{thm}{Theorem}
\newtheorem{KeyLemma}{Key Lemma}
\newtheorem*{theo}{Theorem}

\newtheorem{cor}{Corollary}
\newtheorem{lemma}{Lemma}
\newenvironment{defin}{\medskip\noindent{\sc
Definition}.}{\goodbreak\medskip}
\newenvironment{nota}{\medskip\noindent{\sc
Notations}.}{\goodbreak\medskip}
\newenvironment{remk}{\noindent{\sc
Remark}.}{\goodbreak\vskip10pt}

\newtheorem{ques}{Question}

\def\demo{\medskip\goodbreak\noindent
     \hbox{\sc Proof \kern .3em}\ignorespaces}%
  \def \qedbox{$\square$}%
  \def \qed{\hglue1mm\hfill{\ifmmode\qedbox
     \else\unskip\ \hglue0mm\hfill\qedbox\medskip
      \goodbreak\fi}}%

\def\qed{\hglue1mm\hfill\raise -2pt\hbox{\vrule\vbox to 10pt{\hrule width
4pt
                  \vfill\hrule}\vrule}}

\newcommand{\T}{\mathbb {T}}
\newcommand{\A}{\mathbb {A}}

\newcommand{\R}{\mathbb {R}}

\newcommand{\Z}{\mathbb {Z}}
\newcommand{\N}{\mathbb {N}}

\newcommand{\Cc}{\mathcal {C}}
\newcommand{\Ic}{\mathcal {I}}
\newcommand{\Rc}{\mathcal {R}}
\newcommand{\Mc}{\mathcal {M}}

\newcommand{\Gc}{\mathcal {G}}

\newcommand{\Fc}{\mathcal {F}}

\newcommand{\Sc}{\mathcal {S}}

\begin{document}
\title{The non-hyperbolicity of irrational invariant curves for twist maps and
all that follows}
\author{M.-C. Arnaud\thanks{ANR DynNonHyp ANR BLAN08-2-313375}
\thanks{Avignon University , Laboratoire de Math\'ematiques d'Avignon (EA 2151),  F-84 018 Avignon,
France. e-mail: Marie-Claude.Arnaud@univ-avignon.fr}
\thanks{Member of the Institut Universitaire de France}, P. Berger\thanks{Laboratoire Analyse, G\'eom\'etrie \& Applications
UMR 7539 Institut Galil\'ee Universit\'e Paris 13}
}
\maketitle

\abstract{\em  The key result of this article is\\
{{\textsc{key lemma}:}  if a Jordan curve $\gamma$  is invariant by a given $C^{1+\alpha}$-diffeomorphism $f$ of a surface and if $\gamma$ carries an ergodic hyperbolic probability $\mu$, then  $\mu$ is supported on a periodic orbit.}

From this Lemma we deduce three new results for the $C^{1+\alpha}$ symplectic twist maps $f$ of the annulus:
\begin{enumerate}
\item if $\gamma$ is  a loop at the boundary of an  instability zone such that $f_{|\gamma}$ has  an irrational rotation number, then the convergence of any orbit to $\gamma$  is slower than exponential;
\item if $\mu$ is an invariant probability that is supported in an invariant  curve $\gamma$ with an irrational rotation number, then $\gamma$ is $C^1$ $\mu$-almost everywhere;
\item we prove a part of the so-called ``Greene criterion'', introduced by  J.~M.~Greene in \cite{Gre} in 1978 and never proved:\\
assume that $(\frac{p_n}{q_n})$ is a sequence of rational numbers converging to an irrational number $\omega$; let $(f^k(x_n))_{1\le k\le q_n}$ be a minimizing periodic orbit with rotation number $\frac{p_n}{q_n}$ and let us denote by $\Rc_n$ its mean residue $\Rc_n=\left|1/2-{\rm Tr}(Df^{q_n}(x_n))/4\right|^\frac{1}{q_n}$. Then, if $\displaystyle{\limsup_{n\rightarrow +\infty} \Rc_n>1}$, the   Aubry-Mather set with rotation number $\omega$   is not supported in an invariant curve.\end{enumerate}

} 

\tableofcontents
   \section{Introduction}

    A reason for studying the  positive symplectic twist maps\footnote{all the definitions are given in subsection \ref{ss12}} (PSTM) of the two-dimensional annulus $\A=\T\times\R$  is that they represent (via a
symplectic change of coordinates) the dynamics of a generic symplectic diffeomorphism  of a
surface  near its elliptic periodic points (see for example \cite{Ch1}). One motivating  example of such  a map was introduced by Poincar\'e for the study of  the restricted 3-body problem.\\

The   study of the PSTM was initiated by G.D.~Birkhoff in the 1920s (see \cite{Bir1}). Among other beautiful results, he proved that any essential curve invariant by a symplectic twist map of the annulus is the graph of a Lipschitz map (an {\em essential} curve is a simple loop that is not homotopic to a point ). 

Later, in the '50s, the K.A.M.~theorems provided the existence of some invariant curves  for sufficiently regular symplectic diffeomorphisms of surfaces near their elliptic fixed points (see \cite{Ko}, \cite{Arno}, \cite{Mo} and \cite{Ru}). These theorems provide also some essential invariant curves for the symplectic twist maps that are close to the completely integrable ones. These K.A.M. curves are all very regular (at least $C^3$, see \cite{He2}).  

But general invariant curves for general PSTM have no reason to be so regular. The example of the simple pendulum (see \cite{Ch2}) shows us that an invariant curve can be non-differentiable at one point: the separatrix of the simple pendulum has an angle at the hyperbolic fixed point. In \cite{He2} and \cite{Arna1}, some other examples are given of symplectic twist maps that have a non-differentiable essential invariant curve that contains some periodic points. Some examples of symplectic twist maps  that have an essential invariant curve that is not $C^1$ and that has an irrational rotation number are built in \cite{Arna2} and then improved in \cite{Arna3}. In these examples, the dynamics restricted to the curve is conjugated to the one of a Denjoy counter-example, and the set  where the curve is non differentiable is the orbit of a wandering point (in these examples, the first PSTM is $C^1$ and the second one is $C^2$). 

An open question is then

 \begin{ques}\label{Qgraal} Does there exist an essential non-$C^1$ curve $\gamma$ that is invariant by a PSTM $f$ and such that $f_{|\gamma}$ is minimal?
\end{ques}

Moreover, in the Denjoy examples that we mentionned above, the invariant curve is differentiable along the support of the unique invariant measure supported in the curve. Ke Zhang asked us the following question:
\begin{ques}\label{QKe}
 Is it possible to have some  points of non-differentiability in the support of the invariant measure?
 \end{ques}
 
We will explain  in subsection \ref{ss21} what we mean by ``$C^1$ at one point''. The following result was proved in \cite{Arna1} 
 \begin{theo} {\rm (Arnaud, \cite{Arna1})}
Let $f~:\A\rightarrow \A$ be a $C^1$ PSTM and
let $\eta~:
\T\rightarrow
\R$  be a Lipschitz map the graph of which is invariant by $f$. Let us denote the set of points of $\T$ where $\eta$ is $C^1$ by $U$. Then $U$ contains a dense $G_\delta$ subset of $\T$ that has full Lebesgue  measure.
\end{theo}

Here we give another partial answer to the previous questions (the set where the curve is not $C^1$ is small in a new sense) and  prove

\begin{thm}\label{thregul}
Let $f~:\A\rightarrow \A$ be a $C^{1+\alpha}$ PSTM and
let $\eta~:
\T\rightarrow
\R$  be a Lipschitz map the graph of which is invariant by $f$ such that the rotation number of $f_{|{\rm graph}(\eta)}$ is irrational. Then if $\mu$ is the unique invariant Borel probability measure supported in ${\rm graph}(\eta)$,   the set ${\rm graph}(\eta)$ is $C^1$ on a   subset   that has full $\mu$ measure.
\end{thm} 

The means to obtain Theorem \ref{thregul} is the following Lemma, that is the key ingredient of this article and from which we deduce all the other results. Lemma \ref{thWH} answers in particular to a question that was raised in \cite{Arna1} concerning the existence of Jordan invariant curves that carry a hyperbolic measure. We recall that a Jordan curve is the image of $\T$ by a continuous and injective map.

\begin{KeyLemma}\label{thWH} Let $\gamma:\T\rightarrow \Sc$ be a Jordan curve that is invariant by a $C^{1+\alpha}$ diffeomorphism $f$ of a surface $\Sc$. Assume that $\mu$ is an ergodic Borel probability for $f$  supported in $\gamma$ that is hyperbolic. Then $\mu$ is supported by a periodic orbit and $\gamma$ contains two (stable or unstable) branches of the periodic point.

\end{KeyLemma}
\begin{cor}\label{corWH}
Let $f: \A\rightarrow \A$ be a $C^{1+\alpha}$ PSTM, let $\gamma\subset\A$ be an essential irrational invariant curve for $f$ and let $\mu$ be the unique invariant measure supported in $\gamma$.  Then the Lyapunov exponents of $\gamma$ are zero.\end{cor}

Corollary \ref{corWH} is fundamental to obtain the rate of convergence of the orbits to the boundaries of the bounded zones of instability. Let us recall that when we remove all the essential invariant curves of a given PSTM   $f:\A\rightarrow \A$ from the annulus, we obtain an invariant open set $U$. Following \cite{Arna3}, we call the annular connected components of $U$  the {\em instability zones} of $f$. In \cite{Bir1}, G.~D.~Birkhoff proved that if $U$ is an instability zone, if $U_1$ is a neighborhood of one of its ends (i.e, eventually after compactification,  a connected component of its boundary) and $U_2$ is a neighborhood of the other end, then there exists an orbit  traveling from $U_1$ to $U_2$. This theorem was improved in \cite{Mat2} by J.~Mather who proved that if $\Cc_1$, $\Cc_2$ are the ends of $U$, there exists an orbit whose $\alpha$-limit set is in $\Cc_1$ and $\omega$-limit set is in $\Cc_2$. J.~N.~Mather used variational arguments and after that, P.~Le~Calvez gave in \cite{LeC1} a purely topological proof of this result.

We deduce from Corollary \ref{corWH} and from some results due to A. Furman (see \cite{Fur}) concerning the uniquely ergodic measures the following theorem.

\begin{thm}\label{thrate}
Let $\gamma:\T\rightarrow \A$ be a   loop that is invariant by a $C^{1+\alpha}$ PSTM $f:\A\rightarrow \A$, that has an irrational rotation number and that is at the boundary of an instability zone.
Let $x\in W^s(\gamma)\backslash \gamma$ be a point such that $\displaystyle{\lim_{n\rightarrow +\infty} d(f^n(x), \gamma)=0}$. Then, for all $\varepsilon>0$, we have
$$\lim_{n\rightarrow +\infty} e^{\varepsilon n}d(f^n(x), \gamma)=+\infty.$$
\end{thm} 

In \cite{He2}, M.~Herman proved that $C^r$-generically for any $1\leq r\leq \infty$,  a symplectic twist map has no invariant curve with rational rotation number; he also proved the existence of a non-empty open set of symplectic twist maps having a bounded instability zone. Hence, our result describes what happens in the general case: the boundary curve has an irrational rotation number and all the orbits that converge to this curve converge more slowly than every exponential rate.

Having an upper bound for the rate of convergence to the boundary of an instability zone with an irrational rotation number, we can ask if there exists some lower bound. The only answer that we can give is the answer for the examples that were built in \cite{Arna3}: in these examples, $\gamma$ is a curve that is at the boundary of an instability zone and such that the restricted dynamics $f_{|\gamma}$ is Denjoy, and some orbits are provided such that:
$$\forall n\in \N,  d(f^nx, \gamma)=\frac{d(x,\gamma)}{n(\log n)^{1+\delta}}$$ for $\delta>0$.

Of course, there exist too some (non-generic) examples of boundaries of instability zone that contain      a hyperbolic periodic point and for which an exponential convergence to the boundary may happen. Such an example is given by Birkhoff in \cite{Bir2}.

\begin{cor}\label{corrate}
There exists a $G_\delta$ subset $\Gc$ of a non-empty open set $U$  of PSTM such that every $f\in\Gc$ has an invariant curve $\gamma$ with a non-trivial stable set $W^s(\gamma)\supsetneq \gamma$ such  that the orbit of any point of $W^s(\gamma)\backslash \gamma$ converges to $\gamma$ more slowly than any exponential sequence. 
\end{cor}

\begin{remk}
The phenomenon that is described in Corollary \ref{corrate} points out a new typical dynamical behavior (at least in  $C^\infty$-topology): existence of a large set of dynamics exhibiting a ``slow'' convergence (i.e. non-exponential). In some way, this is reminiscent of Arnol'd diffusion.
\end{remk}

In the 80's, the Aubry-Mather sets were discovered simultaneously and independently by Aubry \& Le Daeron (in \cite{ALD}) and Mather (in \cite{Mat1}). See subsection \ref{ssAM} for details on this now classical theory. These sets are invariant and compact. They are not necessarily on an invariant curve but lie on a Lipschitz graph. We can define for each of these sets a {\em rotation number} and for every real number, there exists at least one   Aubry-Mather set with this rotation number. An Aubry-Mather set with an irrational rotation number can be
 \begin{enumerate}
 \item an invariant curve; 
 \item  the union of a Cantor set and some homoclinic orbits to this Cantor set. Then this Cantor set is not contained in an invariant curve.  \end{enumerate}
 
A  fundamental problem is to identify for which irrational rotation number there exists an invariant curve with this rotation number. A method was proposed by J.~M.~Greene in \cite{Gre}, that was numerically  tested for  the standard map. Let us explain how it works. Let $\rho\in\R$ be an irrational number. We consider a sequence $(\frac{p_n}{q_n})$ of rational numbers that converge to $\rho$ and for each $n$, a minimizing periodic point $x_n=(\theta_n, r_n)$ with rotation number $\frac{p_n}{q_n}$. Then at every $x_n$, the {\em residue} is $r_n=\frac{1}{4}(2-{\rm Tr}(Df^{q_n}(x_n)))$ As noticed by J.~M.~Greene, one advandage of the residue is that it is as regular as $Df$ is and is easily computable (contrarily to the eigenvalues of $Df^{q_n}$).

We are interested in the version of Greene's residue criterion that is presented in \cite{Mac1}. In \cite{Mac1}, R.~S.~MacKay introduces the {\em mean residue} $\Rc_n=|r_n|^\frac{1}{q_n}$ and gives the following conjecture\\
{\textsc{Residue criterion}.} {\em Let $(\frac{p_n}{q_n})$ be a good sequence of rational numbers that converges to $\rho$ and let $x_n$ be a minimizing periodic orbit with rotation number $\frac{p_n}{q_n}$ and mean residue $\Rc_n$. Then $\displaystyle{\lim_{n\rightarrow \infty}\Rc_n=\mu(\rho)}$ exists. If $\mu(\rho)\leq 1$, there exist an invariant curve with rotation number $\rho$ and if $\mu(\rho)>1$, such a curve does not exist.
}\\
Some partial results are proved in \cite{Mac1}. Using similar ideas and key Lemma \ref{thWH}, we will prove

\begin{thm}\label{thgreene}
Let $f:\A\rightarrow \A$ be a PSTM and let $\rho$ be an irrational number. Let $(x_n)$ be a sequence of minimizing periodic points with rotation number $\frac{p_n}{q_n}$ and mean residue $\Rc_n$ such that $\displaystyle{\limsup_{n\rightarrow \infty}\Rc_n>1}$. Then there exists no essential invariant curve with rotation  number $\rho$. Moreover, the  Aubry-Mather set $K$ with rotation number $\rho$ contains  a Cantor $C$ set that carries a hyperbolic invariant probability  and that is $C^1$-irregular $\mu$-almost everywhere.
\end{thm}

In this case, a subsequence of the sequence of  periodic orbits  $(\{f^i(x_n): \; i\ge 0\})_n$ converges  for the Haurdorff metric to an invariant  subset of $K$ that contains the Cantor set $C$.   We proved in \cite{Arna1} that because $C$ carries an invariant hyperbolic probability measure denoted by $\mu$, the compact set $C$ is $C^1$-irregular $\mu$-almost everywhere. A natural question is then.  

 \begin{ques} Can we see the $C^1$-irregularity of $C$ with the help of the orbits of the $x_n$ (with a computer)?\end{ques}

\subsection{Structure of the article}
After explaining what are the Lyapunov charts, we will prove key Lemma  \ref{thWH} and Corollary \ref{corWH}  in section \ref{section2}.\\
 In section \ref{section3}, we state Theorem \ref{thzerolyap} about the rate of attraction of the uniquely ergodic measures. Joined with  key Lemma  \ref{thWH}, Theorem \ref{thzerolyap}  implies Theorem \ref{thrate} and Corollary \ref{corrate}.\\
 In section \ref{s4}, we will recall some facts about the Green bundles and then  we will prove Theorem \ref{thregul}.\\
In section \ref{s5}, we will prove Theorem \ref{thgreene}.
\subsection{Some notations and definitions}\label{ss12}
Before giving the proofs, let us introduce some notations and definitions.

\begin{nota}
\noindent $\bullet$ $\T=\R/\Z$ is the circle.

\noindent $\bullet$ $\A=\T\times \R$ is the annulus and an element of $\A$ is denoted
by $x=(\theta, r)$.

\noindent $\bullet$ $\A$ is endowed with its usual symplectic form, $\omega=d\theta\wedge dr$ and its usual Riemannian metric.

\noindent  $\bullet$ $\pi: \T \times \R \rightarrow\T$ is the first projection and its lift is denoted by $\Pi: \R^2\rightarrow \R$.

\noindent $\bullet$ For every $x\in \A$, $V(x)=\ker D\pi(x)\subset T_x\A$ is the linear vertical at $x$.

\noindent$\bullet$ $p:\R^2\rightarrow \A$ is the universal covering.
\end{nota}

\begin{defin} A $C^1$ diffeomorphism $f: \A\rightarrow \A$ of the annulus that is isotopic
to identity  is a {\em positive twist map} (resp. {\em negative twist map}) if there exists $\varepsilon>0$ such that  for any $x\in \mathbb A$, we have: $D(\pi\circ f)(x)(0, 1)>\varepsilon$ (resp. $D(\pi\circ f)(x)(0, 1)<-\varepsilon$).
 A {\em
twist map} may be positive or negative. 
\end{defin}
We recalled in the introduction that any essential curve that is invariant by a symplectic twist map is the graph of a Lipschitz map. 

\begin{defin}
Let $\Gamma$ be an essential invariant curve of a symplectic twist map $f$ of the annulus. Then,  if we project the restricted dynamics to $\Gamma$ on the circle, we obtain an orientation preserving homeomorphism of the circle, that has a rotation number. If this rotation number is irrational, we will say that the curve is {\em irrational}. 
\end{defin}
Let us recall that the dynamics restricted to an irrational invariant curve is uniquely ergodic.


\begin{defin}
Let $f:\A\rightarrow \A$ be a $C^1$-diffeomorphism. Let $\mu$ be a Borel probability that is left invariant by $f$. We say that $\mu$ is a {\em hyperbolic measure} if it has one negative Lyapunov exponent and one positive Lyapunov exponent.
\end{defin}

\begin{defin}
Let $K\subset \A$ be a non-empty compact subset that is invariant by a diffeomorphism $f:\A\rightarrow \A$. The {\em stable set} of $K$ for $f$ is
$$W^s(K, f)=\{ x\in\A; \lim_{n\rightarrow +\infty}d(f^n(x), K)=0\}.$$
\end{defin}

\begin{defin}
\begin{enumerate}
\item Let $K\subset \A$ be a non-empty compact subset  and let $x\in K$ be a point of $K$. The {\em (Bouligand) paratangent cone} to $K$ at $x$ if the subset $P_xK$ of $T_x\A$ that contains all the limit points of the sequences
$$t_n(x_n-y_n)$$
where $t_n\in\R$ and $(x_n)$, $(y_n)$ are two sequences of points of $K$ that converge to $x$.
\item The set $K$ is {\em $C^1$-regular} (in fact we will say $C^1$) at $x$ if $P_xK$ is contained in a line.
\end{enumerate}
\end{defin}
\begin{remk}
\begin{enumerate}
\item A loop in $\A$ is $C^1$-regular if and only if it is $C^1$ as a submanifold.
\item The graph of a Lipschitz map $\eta:\T\rightarrow \R$ is $C^1$-regular if and only if $\eta$ is $C^1$.
\end{enumerate}
\end{remk}

\subsection{Aubry-Mather theory}\label{ssAM}
For proofs, see \cite{Ban} and  \cite{Gol}. In this section, we assume that the twist maps that we consider are a little more than symplectic: they are {\em exact symplectic } i.e. $f^*(rd\theta)-rd\theta$ is exact.\\
To any lift $F:\R^2\rightarrow \R^2$ of a positive exact symplectic twist map (PSTM) $f:\A\rightarrow \A$	we may associate a $C^2$  {\em generating function} $S:\R^2\rightarrow \R$ that satisfies  
\begin{enumerate}
\item[$\bullet$] $S(\theta+1, \Theta+1)=S(\theta,  \Theta)$;
 \item[$\bullet$] there exists $\varepsilon>0$ such that: $\frac{\partial^2 S}{\partial \theta\partial \Theta}(\theta,  \Theta) \leq -\varepsilon$;
 \item[$\bullet$] $F$ is implicitly given by:
 $$F(\theta,r)=(\Theta, R)\Longleftrightarrow \left\{\begin{matrix}r=-\frac{\partial S}{\partial \theta}(\theta,  \Theta)\\ R=\frac{\partial S}{\partial \Theta}(\theta,  \Theta)
 \end{matrix}\right. .$$ 
 \end{enumerate}
For every $k\geq 2$, $\theta_0, \theta_k\in\R$, the function $\Fc_{(0, \theta_0), (k, \theta_k)}=\Fc:\R^{k-1}\rightarrow \R$ is defined by\\
 $\displaystyle{\Fc(\theta_1, \dots, \theta_{k-1})=\sum_{j=1}^kS(\theta_{j-1}, \theta_j)}$. The function $\Fc_{(0, \theta_0), (k, \theta_k)}$ has a minimum, and at every critical point for $\Fc$, the following sequence is a piece of orbit for $F$:
 $$(\theta_0, -\frac{\partial S}{\partial \theta}(\theta_0, \theta_1)), (\theta_1,  \frac{\partial S}{\partial \Theta}(\theta_0, \theta_1)), (\theta_2,    \frac{\partial S}{\partial \Theta}(\theta_1, \theta_2)), \dots , (\theta_k,  \frac{\partial S}{\partial \Theta}(\theta_{k-1}, \theta_k)).$$
 \begin{defin} An orbit $(\theta_n, r_n)$ of $F$ (and by extension its projection on $\A$) is minimizing if every finite segment $(\theta_n)_{\ell\leq n\leq k}$ is {\em }minimizing for $\Fc_{(\ell-1, \theta_{\ell-1}), (k+1, \theta_{k+1})}$.\end{defin}
 It can be proved that every minimizing orbit $(\theta_n, r_n)$ has a rotation number $$\displaystyle{\rho=\lim_{n\rightarrow \pm \infty} \frac{\theta_n-\theta_0}{n}}.$$
 
The set of points $(\theta, r)\in\R^2$ having a minimizing orbit is denoted by $\Mc$. Then it is closed and its projection $p(\Mc)\subset \A$ is closed too. The rotation number $\rho:\Mc\rightarrow \R$ is continuous and for every $\alpha\in\R$, the set $\Mc_\alpha=\{ x\in\Mc, \rho(x)=\alpha\}$ is non-empty.\\

If $\alpha$ is irrational, then $K_\alpha=p(\Mc_\alpha)\subset \A$ is the graph of a Lipschitz map above a compact subset of $\T$. Moreover, there exists a bilipschitz orientation preserving homeomorphims $h:\T\rightarrow \T$ such that 
$$\forall x\in K_\alpha, h(\pi(x))=\pi(f(x)).$$
hence $K_\alpha$ is:
\begin{enumerate}
\item[-]  either not contained in an invariant loop and then is the union of a Cantor set $C_\alpha$ on which the dynamics is minimal and some homoclinic orbits to $C_\alpha$;
\item[-] or is an invariant loop. In this case the dynamics restricted to $K_\alpha$ can be minimal or Denjoy.
\end{enumerate}
If $\alpha=\frac{p}{q}$ is rational, then $\Mc_\alpha$ is the disjoint union of 3 sets:
\begin{enumerate}
\item[$\bullet$] $\Mc_\alpha^{\rm per}=\{ x\in \Mc_\alpha, \Pi\circ F^q(x)=\Pi(x)+p\}$;
\item[$\bullet$] $\Mc_\alpha^+=\{ x\in \Mc_\alpha, \Pi\circ F^q(x)>\Pi(x)+p\}$;
\item[$\bullet$] $\Mc_\alpha^-=\{ x\in \Mc_\alpha, \Pi\circ F^q(x)<\Pi(x)+p\}$.
\end{enumerate}
The two sets $K_\alpha^+=p(\Mc_\alpha^{\rm per}\cup\Mc_\alpha^+)$ and $K_\alpha^-=p(\Mc_\alpha^{\rm per}\cup\Mc_\alpha^-)$ are then invariant  Lipschitz graphs above a compact part of $\T$. the points of $p(\Mc_\alpha^+\cup \Mc_\alpha^-)$ are   heteroclinic orbits to some peridoc points contained in $p(\Mc_\alpha^{\rm per})$.

\begin{defin}
If $\alpha$ is irrational, the {\em Aubry-Mather set} with rotation number $\alpha$ is $\Mc_\alpha$;.\\
If $\alpha$ is rational, the two {\em Aubry-Mather sets} with rotation number $\alpha$ are $K_\alpha^-$ and $K_\alpha^+$. 
\end{defin}

\begin{remk} If $x$ is a periodic point of a PSTM $f$ that is contained in a Aubry-Mather set, then there exists a line bundle $G_-\subset T\A$ along the orbit of $x$ that is transverse to the vertical fiber and   invariant by $Df$.  This line bundle is described in subsection \ref{ssGR} and is one of the two Green bundles. Moreover,  the restricted dynamics $Df_{|G_-}$ is orientation preserving. Indeed, it is proved in the proof of proposition 7 in \cite{Arna5} that $D(\pi\circ f)_{|G_-}=b\Delta s$ where $Df=\begin{pmatrix}a&b\\ c&d\end{pmatrix}$ (we have $b>0$ because $f$ is a PSTM) and $\Delta s =s_--s_{-1}$ where $s_-$ is the slope of $G_-$ and $s_{-1}$ is slope of the inverse image of the vertical $Df^{-1}.V$ (as noticed in \cite{Arna5}, we have $s_->s_{-1}$). Hence if $\tau$ is the period of $x$, $Df^\tau(x)$ has a positive eigenvalue.
\end{remk}

\section{On the hyperbolic measures that are supported in an essential invariant curve}\label{section2}
 Let $\gamma:\T\rightarrow \Sc$ be a Jordan curve of a surface $\Sc$ that is invariant by a $C^{1+\alpha}$ diffeomorphism $f:\Sc\rightarrow \Sc$. Assume that $\mu$ is an ergodic   Borel probability for $f$  that is hyperbolic  with support in $\gamma$. \\
We want to prove that  $\mu$ is supported by a periodic orbit and $\gamma$ contains two (stable or unstable) branches of the periodic cycle.

\begin{remk}
If $f_{|\gamma}$ is not orientation preserving, we replace $f$ by $f^2$ (note that then $\mu$ must be supported by a fixed point of $f^2$).
\end{remk}
\subsection{Pesin theory and Lyapunov charts}\label{ss21}
We denote by $d$ a Riemannian distance in $\Sc$, by $\|.\|$ the associated norm  in $T\Sc$ and by $|.|$ the norm sup defined on $\R^2$ by $| (x, y)|=\sup\{ |x|, |y|\}$.\\
We recall the terminology and the results on Pesin theory that are given in \cite{L-Y} (see \cite{K-H} too).\\
Let $ \Gamma$ be the set of  regular points for $\mu$, i.e. the subset of points $x\in {\rm supp}\mu$  where there exists a splitting $T_x\Sc=E^s(x)\oplus E^u(x)$ such that
$$\forall v\in E^s(x)\backslash\{0\}, \lim_{k\rightarrow \pm\infty}\frac{1}{k}\log\left( \| Df^k(x)v\|\right)=-\lambda_1$$\quad{\rm and}\quad $$\forall v\in E^u(x)\backslash\{0\}, \lim_{k\rightarrow \pm\infty}\frac{1}{k}\log\left( \| Df^k(x)v\|\right)=\lambda_2;$$
where $\lambda_1$, $\lambda_2$ are positive real numbers. We introduce the notation $\lambda=\min\{\lambda_1, \lambda_2\}$.\\
For any $\rho>0$, we denote the square $[-\rho, \rho]^2$ by $R(\rho)$. Let $\varepsilon\in (0, \frac{\lambda}{10})$ be given.\\
 We can define a local chart in some neighborhood of any regular point; the size of the neighborhood, the local chart and the size of the estimates varies with $x\in\Gamma$ in a measurable way. \\
 More precisely, there is a measurable function $\ell:\Gamma\rightarrow [1, \infty)$ such that $e^{-\varepsilon}\ell(x)\leq\ell(f(x))\leq e^\varepsilon \ell(x)$, a constant $K>0$  and a $C^\infty$ embedding $\Phi_x:R(\frac{1}{\ell(x)})\rightarrow \Sc$ with the following properties:
 \begin{enumerate}
 \item[(i)] $\Phi_x(0)=x$; $D\Phi_x(0)(\R\times \{ 0\})=E^u(x)$ and $D\Phi_x(0)(\{0\}\times \R)=E^s(x)$;
 \item[(ii)] if we denote by $F_x=\Phi_{f(x)}\circ f\circ \Phi_{x}^{-1}$ the connecting map between the chart at $x$ and the chart at $f(x)$ that are defined whenever it makes sense and similarly $F_x^{-1}=\Phi_{f^{-1}x}^{-1}\circ f^{-1}\circ \Phi_{x}$, then we have:
 $$\left | DF_x(0)\begin{pmatrix} 1\\ 0\end{pmatrix}\right |\geq e^{\lambda_2-\varepsilon }\left | \begin{pmatrix} 1\\ 0\end{pmatrix}\right|=e^{\lambda_2-\varepsilon }$$
 and 
 $$\left | DF_x(0)\begin{pmatrix} 0\\ 1\end{pmatrix}\right |\leq e^{-(\lambda_1-\varepsilon) }\left | \begin{pmatrix} 0\\ 1\end{pmatrix}\right |=e^{-(\lambda_1-\varepsilon )}$$
 \item[(iii)] if $L(g)$ denotes the Lipschitz constant of a function $g$, then
 $$L(F_x-DF_x(0))\leq \varepsilon, L(F_x^{-1}-DF_x^{-1}(0))\leq\varepsilon\quad{\rm and}\quad L(DF_x), L(DF_x^{-1})\leq \ell(x);$$
 \item[(iv)] for all $z, z'\in R(\frac{1}{\ell(x)})$, we have
 $$\frac{1}{K}d(\Phi_x(z), \Phi_x(z'))\leq |z-z'|\leq \ell(x)d(\Phi_x(z), \Phi_x(z')).$$

 \end{enumerate}
 We deduce from  (ii) and (iii) that there exists $\Lambda>0$ such that for all $x\in \Gamma$, we have
 $$F_x\left( R(\frac{e^{-(\Lambda+\varepsilon)}}{\ell(x)})\right)\subset R(\frac{1}{\ell(fx)})\quad{\rm and}\quad F_x^{-1}\left( R(\frac{e^{-(\Lambda+\varepsilon)}}{\ell(x)})\right)\subset R(\frac{1}{\ell(f^{-1}x)})$$
 From now we will use the small charts $\left( R(\frac{e^{-(\Lambda+\varepsilon)}}{\ell(x)}), \Phi_x\right)$, because $F_x$ and $F_x^{-1}$ are defined on the whole domain of these charts, but to be short we will use the notation $\ell(x)$ instead of $e^{\Lambda+\varepsilon}\ell(x)$. We will call them (Lyapunov) $(\varepsilon, \ell)$-charts.\\
 
 We recall that the stable (resp. unstable) manifold of $x\in \Gamma$ is defined by:
 $$W^s(x)=\{ y\in \Sc; \limsup_{n\rightarrow +\infty} \frac{1}{n}\log d(f^nx, f^ny)<0\}$$
(resp. 
 $$W^u(x)=\{ y\in \Sc; \limsup_{n\rightarrow +\infty} \frac{1}{n}\log d(f^{-n}x, f^{-n}y)<0\}).$$
 The local stable (resp. unstable) manifold at $x$ associated with $(\Phi_x)$ is then the connected component of $W^s(x)\cap \Phi_x(R(\frac{1}{\ell(x)}))$ (resp. $W^u(x)\cap \Phi_x(R(\frac{1}{\ell(x)}))$) that contains $x$ and  is denoted by $W^s_{\rm loc}(x)$ (resp. $W^u_{\rm loc}(x)$) . The $\Phi_x^{-1}$ images of these sets are denoted by $W^s_x(0)$ and $W_x^u(0)$.\\
 The function $\ell$ being eventually replaced by $\delta\ell$ for some large $\delta>0$, it can be proved that $W^u_x(0)=\{ (t, g_x^u(t)); t\in [-\frac{1}{\ell(x)}, \frac{1}{\ell(x)}]\}$ where $g_x^u:[-\frac{1}{\ell(x)}, \frac{1}{\ell(x)}]\rightarrow [-\frac{1}{\ell(x)}, \frac{1}{\ell(x)}]$ is a $C^{1+\alpha}$ function such that $g_x^u(0)=0$ and $|g_x^{u\prime}(0)|\leq \frac{1}{3}$ and that $W^s_x(0)=\{ (  g_x^s(t),t); t\in [-\frac{1}{\ell(x)}, \frac{1}{\ell(x)}]\}$ where $g_x^s:[-\frac{1}{\ell(x)}, \frac{1}{\ell(x)}]\rightarrow [-\frac{1}{\ell(x)}, \frac{1}{\ell(x)}]$ is a $C^{1+\alpha}$ function such that $g_x^s(0)=0$ and $|g_x^{s\prime}(0)|\leq \frac{1}{3}$.

\subsection{Proof of key Lemma  \ref{thWH}}
\begin{figure}[!h]
    \centering
    \includegraphics{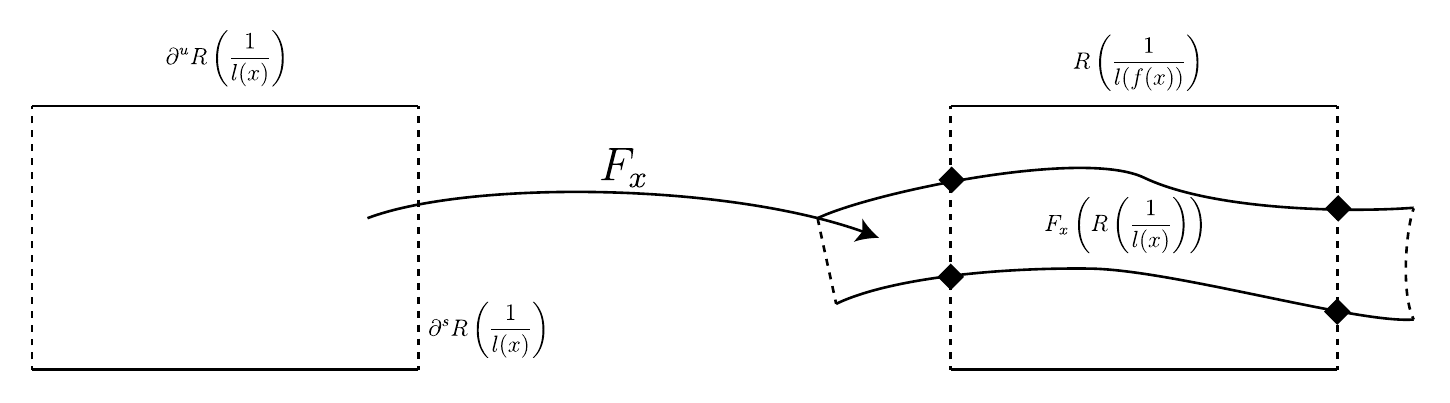}
    \label{figure}
  \end{figure}

We eventually change $\ell$ into $\delta\ell$ for a large $\delta$ to be sure that $\gamma$ is contained in no $\Phi_x(R(\frac{1}{\ell (x)}))$.\\
We decompose the boundary $\partial R(\rho)$ of $R(\rho)$ into $\partial^s R(\rho)=\{-\rho, \rho\}\times [-\rho, \rho]$ and $\partial^u R(\rho)=[-\rho, \rho]\times \{ -\rho, \rho\}$ (see Figure \ref{figure}). \\

Then $F_x(R(\frac{1}{\ell(x)}))$ is a curved square and we have: $F_x(R (\frac{1}{\ell(x)}))\cap \partial ^u R(\frac{1}{\ell(fx)})=\emptyset$ and $F_x(\partial ^u R(\frac{1}{\ell(x)}))\cap \partial R(\frac{1}{\ell(fx)})$ contains 4 points, two on the right component of $\partial^s R(\frac{1}{\ell(x)})$ and two on its left component.

The loop $\gamma$ is endowed with an orientation. We shall assume, up to considering $f^2$ instead of $f$, that $f|\gamma$ preserves this orientation.
For every $x\in \Gamma$, we denote by $\gamma_x$ the component of $\gamma\cap \Phi_x(R(\frac{1}{\ell(x)}))$ that contains $x$. 
Then $\gamma_x\backslash \{ x\}$ has two connected components, and we denote by $\eta_x$ the closure of the one after $x$ following the orientation (the same reasoning remains the same for the component before $x$). We observe that $f(\eta(x))\cap \eta(f(x))\supsetneq \{x\}$. 
We will prove that $x$ is a periodic point and that $\eta_x\subset W^s_{\rm loc}(x)$ or $\eta_x\subset W^u_{\rm loc}(x)$. This will give the conclusion of key Lemma  \ref{thWH}. \\
The set  $C_x=\phi_{x}^{-1}(\eta_x)$ is an arc (i.e. the image of $[0, 1]$ by a (continuous) embedding) that joins $0=C_x(0)$ to $C_x(1)\in\partial R(\frac{1}{\ell(x)})$. There are two cases: $C_x(1)\in \partial^s R(\frac{1}{\ell(x)})$ or $C_x(1)\in \partial ^u R(\frac{1}{\ell(x)})$. 
\begin{lemma}
We have either for $\mu$ almost every $x\in \Gamma$, $C_x(1)\in \partial^s R(\frac{1}{\ell(x)})$ or for $\mu$ almost every $x\in \Gamma$, $C_x(1)\notin \partial^s R(\frac{1}{\ell(x)})$ (and then $C_x(1)\in \partial ^u R(\frac{1}{\ell(x)})$). 
\end{lemma}
\demo
Let us assume for an $x\in\Gamma$  that $C_x(1)\in \partial^s R(\frac{1}{\ell(x)})$. Then $F_x\circ C_x$ joins $0$ to $F_x(C_x(1))\in F_x(R(\frac{1}{\ell(x)}))\backslash R(\frac{1}{\ell(fx)})$. This implies that $C_{fx}\subset F_x(C_x)$ joins $0$ to a point of $C_{fx}(1)\in\partial^s R(\frac{1}{\ell(fx)})$. In  a similar way, we obtain: $\forall n\geq 0, C_{f^nx}(1)\in \partial^s R(\frac{1}{\ell (f^nx)})$.

The map $\Ic:\Gamma\rightarrow \{ 0, 1\}$ is defined by $\Ic(x)=0$ if $C_x(1)\notin \partial^s R(\frac{1}{\ell (x)})$ and by $\Ic(x)=1$ if $C_x(1)\in \partial^s R(\frac{1}{\ell (x)})$. Then $\Ic$ is measurable and we just proved that $\Ic$ is  non-decreasing along the orbits. Hence $\Ic\circ f-\Ic\leq 0$. As $\int(\Ic\circ f-\Ic)d\mu=0$, we deduce that we have $\mu$-almost everywhere $\Ic\circ f=\Ic$. Because $\mu$ is ergodic, $\Ic$ is constant $\mu$-almost everywhere and this gives the wanted result.


From now we assume that we have almost everywhere $C_x(1)\in \partial^s R(\frac{1}{\ell (x)})$.  Let us recall that $W^u_x(0)$ is the graph of $g_x^u$ and let us use the notation $C_x(t)=(c^1_x(t), c^2_x(t))$.  If $\delta(x)=\max\{ |c^2_x(t)-g_x^u(c^1_x(t))|; t\in [0,1]\}$, we use the (ii) and the (iii) of subsection \ref{ss21} to deduce that $\delta(fx)\leq e^{\lambda_1-3\varepsilon}\delta(x)$. Hence $$\int\delta(x)d\mu(x)= \int\delta(fx)d\mu(x)\leq e^{\lambda_1-3\varepsilon}\int\delta(x)d\mu(x)$$ and thus $\delta=0$ $\mu$-almost everywhere. This implies that the corresponding branch of $W^u(x)$ (and then its orbit) is in $\gamma$. 

Assume that the rotation number of $f_{|\gamma}$ is not rational. Then either the dynamics $f_{|\gamma}$  is minimal and we have never for two different points $x\not=y$ $(*) \displaystyle{\lim_{n\rightarrow +\infty} d(f^{-n}x, f^{-n}y)=0}$   or the dynamics is Denjoy and $(*)$ happens only in the  wandering intervals, i.e. for no $x\in{\rm supp}\mu$ but the one that are endpoints of these wandering intervals. The numbers of wandering intervals being countable, we obtain a contradiction. Hence $x$ has to be periodic.

\section{On the rate of convergence to a boundary of an instability zone with irrational rotation number}\label{section3}
\subsection{A result for the rate of convergence to a uniquely ergodic measure with zero Lyapunov exponents}
Let us assume that $\gamma$ is an essential invariant curve by a PSTM $f:\A\rightarrow \A$ that is at the boundary of an instability zone. We assume too that the rotation number of $f_{|\gamma}$ is irrational. Then $f_{|\gamma}$ is uniquely ergodic and by Corollary  \ref{corWH}, its unique invariant probability has zero Lyapunov exponents. Hence Theorem \ref{thrate} is just a consequence of the following theorem.

\begin{thm}\label{thzerolyap}
Let $f:M\rightarrow M$ be a $C^1$-diffeomorphism of a manifold $M$. Let $K\subset M$ be a compact set that is invariant by $f$. We assume that $f_{|K}$ is uniquely ergodic and we denote the unique Borel invariant probability with support in $K$ by $\mu$. We assume that all the Lyapunov exponents of $\mu$ are zero. Let $x_0\in W^s(K, f)\backslash K$. Then we have:

$$ \forall \varepsilon>0, \lim_{n\rightarrow +\infty} e^{\varepsilon n}d(f^n(x_0), K)=+\infty.$$
\end{thm}
Let us now prove this theorem. 
\begin{proof}By hypothesis, we have for $\mu$-almost every point :
$$\lim_{n\rightarrow \pm\infty} \frac{1}{n}\log\| Df^n(x)\|=0.$$
We can use a refinement Kingman's subadditive ergodic theorem that is due to A.~Furman (see Corollary 2 in \cite{Fur}) that implies that we have 
$$ 
\limsup_{n\rightarrow \pm\infty} \max_{x\in K}\frac{1}{n}\log \| Df^n(x)\|\leq 0
.$$

In particular, for any $\varepsilon>0$, there exists $N\geq 1$ such that:
\begin{equation}\label{E1}\forall x\in K, \forall n\geq N, \frac{1}{n}\log\| Df^{-n}(x))\| \leq \frac{\varepsilon}{8}.\end{equation}

Observe that for the following norm with $k\ge N$ large:
\[\|u\|'_x= \sum_{n=0}^k e^{-n \varepsilon /4}\| Df^{-n}(x)u\|_x,\]
satisfies uniformly on $x$ for $u\not=0$:
\[\frac{\|D f^{-1}(x) u\|'_{f^{-1}(x)}}{\|u\|'_x}= e^{\varepsilon /4} + \frac{e^{-k \varepsilon /4}\| Df^{-k-1}(x)u\|_x-e^{\varepsilon/4}\|u\|}{\|u\|'_x}\]
\[\le e^{\varepsilon /4} + \frac{e^{-k \varepsilon /4}\| Df^{-k-1}( x)u\|_x}{\|u\|'_x}\le e^{\varepsilon /4} +e^{-k\varepsilon /8}\]
Hence by changing the Riemannian metric by the latter  one, we can assume that the norm of $D_xf^{-1}$ is smaller than $e^{ \varepsilon /3}$ for every $x\in K$.

Consequently, on a $\eta$-neighborhood $N_\eta$ of $K$, it holds for every $x\in N_\eta$ that:
\[\|D_x f^{-1}\|'\le e^{\varepsilon /2}\]

Let $x_0\in M$ be such that $x_n:=f^n(x_0)\to K$, we want to show that $$\liminf \frac1n \log d(x_n,K)\ge -\varepsilon .$$  
We suppose that $\liminf \frac1n \log d(x_n,K)< -\varepsilon $ for the sake of a contradiction. 
Hence there exists $n$ arbitrarily large so that $x_n$ belongs to the $e^{-n\varepsilon }\eta$-neighborhood of $K$.
Let $\gamma$ be a $C^1$-curve connecting $x_n$ to $K$ and of length at most $e^{-n\varepsilon }\eta$. 
By induction on $k\le n$, we notice that $f^{-k}(\gamma)$ is a curve that connects $x_{n-k}$ to $K$, and has length at most $e^{-n\varepsilon+k\varepsilon/2 }\eta$, and so is included in $V_\eta$. Thus the point $x_0$ is at most $e^{-n\varepsilon/2 }\eta$-distant from $K$. Taking $n$ large, we obtain that $x_0$ belongs to $K$. A contradiction.

\end{proof}

\subsection{Proof of Corollary \ref{corrate}}
 Let $U$ be an open non empty set of symplectic twist maps where the K.A.M. theorems can be used: there are at least two irrational, essential invariant curves. For $f$ in U, we denote two such curves by $\gamma_1(f)$ and $\gamma_2(f)$ and by $\rho_1(f)$, $\rho_2(f)$ their rotation numbers. 
 
 We consider then the dense $G_\delta$ subset $\Gc$ of $U$ of symplectic twist maps such that every minimizing periodic orbit is hyperbolic and all the heteroclinic intersection points between such minimizing periodic points are transverse. M. Herman proved in \cite{He2} that an element of $\Gc$ has no invariant curve that contains a periodic point. 
 
 Let $f\in\Gc$ and let us consider a rational number $\frac{p}{q}$ that is between $\rho_1(f)$ and $\rho_2(f)$. Let $x$ be a minimizing periodic point with rotation number $\frac{p}{q}$. Then $x$ is in the bounded part of $\A$ that is between $\gamma_1(f)$ and $\gamma_2(f)$. 
 
 As $f\in\Gc$, there is no invariant curve with rotation number $\frac{p}{q}$. As the set of invariant curve is closed (see \cite{Bir1}) and the rotation number depends continuously of the invariant curve, we find two numbers $\omega_1<\frac{p}{q}<\omega_2$ such that $f$ has no curve with rotation number in $]\omega_1, \omega_2[$ but $f$ has an invariant  curve $\eta_1$ with rotation number $\omega_1$ and an invariant curve $\eta_2$ with rotation number $\omega_2$. As $f\in\Gc$, the numbers $\omega_1$ and $\omega_2$ are irrational and they are at the boundary of an instability zone (that contains $x$). Hence we can apply Theorem \ref{thrate} to conclude.
 \section{Size of the set of $C^1$-regularity of an irrational invariant curve by a symplectic twist map}\label{s4}
 \subsection{Green bundles and regularity}\label{ssGR}
 A reference for what is in this subsection is \cite{Arna4} and \cite{Arna1}.
 
 \begin{defin} Let $(\theta, r)$ be a point of minimizing orbit for  a PSTM $f$. Then the two Green bundles at $(\theta, r)$ are defined by:
$$G_+(\theta, r)=\lim_{n\rightarrow +\infty} Df^n\cdot V(f^{-n}(\theta, r))\quad {\rm and}\quad G_-(\theta, r)=\lim_{n\rightarrow +\infty} Df^{-n}\cdot V(f^{n}(\theta, r)).$$
\end{defin}
These two Green bundles are measurable, invariant by the differential $Df$ and transverse to the linear vertical $V$.

\begin{nota} \begin{enumerate}\item We denote by   $s_-$, $s_+$ the  slopes of the two Green bundles: $$G_\pm(\theta, r)=\{(t, ts_\pm (\theta, r)); t\in\R\}\subset \R^2=T_{(\theta, r)}\A.$$
\item  If $K\subset \A$ is a compact subset that is contained in some Lipschitz graph (for example $K$ can be an Aubry-Mather set), if $x\in K$, we denote by $p_xK$ the set of all the slopes of the elements of the paratangent cone $P_xK$ (that was defined in subsection \ref{ss12})
$$P_xK\backslash\{ 0\}=\{ (t, tp); t\in\R^*, p\in p_xK\}.$$
\item if $A, B\subset \R$, the relation $A\leq B$ means $\forall a\in A, \forall b\in B, a\leq b$.
\end{enumerate}
\end{nota}
Note that $K$ is $C^1$ at $x$ if and only if $p_xK$ has one or zero elements.

We proved in \cite{Arna4} that if $K$ is a minimizing Aubry-Mather set for the PSTM $f$, then  we have 
$$\forall (\theta, r)\in  K,\quad  s_-(\theta,r)\leq p_{(\theta, r)}K\leq s_+(\theta,r).$$
Hence, to prove that the minimizing Aubry-Mather  set is $C^1$ at a $(\theta, r)\in K$, we only have to prove (but this is not a equivalence) that $s_-(\theta,r)=s_+(\theta,r)$.\\
We proved in \cite{Arna1} the following result:
\begin{theo} {\bf (Arnaud, \cite{Arna1})}
Let $f~:\A\rightarrow \A$ be a $C^1$ PSTM and
let $\eta~:
\T\rightarrow
\R$  be a Lipschitz map the graph of which graph is invariant by $f$. Then the set $\{ \theta\in \T; s_-(\theta, \eta(\theta))=s_+(\theta, \eta(\theta))\}$ is   a dense $G_\delta$
subset $U$ of $\T$ with Lebesgue measure 1. Hence the map $\eta$ is $C^1$ on a dense $G_\delta$ subset of $\T$ that has full Lebesgue measure.
\end{theo}
\subsection{Green bundles and Lyapunov exponents}\label{ssGrLyap}
We proved in \cite{Arna4} that if an invariant measure $\mu$ is supported in a minimizing Aubry-Mather set then
\begin{enumerate}
\item[$\bullet$] either it is hyperbolic (i.e. its Lyapunov exponents are non-zero) and $\mu(\{ s_-\not=s_+\})=1$;
\item[$\bullet$] or its Lyapunov exponents are zero and  $\mu(\{ s_-=s_+\})=1$.
\end{enumerate}
\subsection{Proof of Theorem \ref{thregul}}

Theorem \ref{thregul} is a consequence of the following theorem.
\begin{thm}\label{thPP}
Let $f~:\A\rightarrow \A$ be a $C^{1+\alpha}$ PSTM  and
let $\eta~:
\T\rightarrow
\R$  be a Lipschitz map the graph of which  is invariant by $f$ such that the rotation number of $f_{|{\rm graph}(\eta)}$ is irrational. Then if $\mu$ is the unique invariant Borel probability measure supported in ${\rm graph}(\eta)$, we have $\mu(\{ \theta\in \T; s_-(\theta)=s_+(\theta)\})=1$.  Hence the map $\eta$ is $C^1$ on a   subset of $\T$ that has full $\mu$ measure.
\end{thm}
Theorem \ref{thPP} is a straightforward consequence of 
 key Lemma  \ref{thWH} and of the results that are contained in subsection \ref{ssGrLyap}.
\section{About Greene's criterion}\label{s5}
Let us assume that $f:\A\rightarrow \A$ is a PSTM, that $\rho$ is an irrational number, that $(x_n)$ is  a sequence of minimazing periodic points with rotations numbers $\frac{p_n}{q_n}$ such that $\displaystyle{\lim_{n\rightarrow +\infty} \frac{p_n}{q_n}=\rho}$. We denote the mean residue of $x_n$ by $\Rc_n=\left|\frac{2-{\rm Tr}(Df^{q_n}(x_n))}{4}\right|^{\frac{1}{q_n}}$ and we assume that $\displaystyle{\lim_{n\rightarrow +\infty}\Rc_n>1}$. 

If $\displaystyle{\mu_n=\frac{1}{q_n}\sum_{j=1}^{q_n}\delta_{f^j(x_n)}}$ is the measure equidistributed along the orbit of $x_n$, the sequence $(\mu_n)$ converges to the unique invariant measure $\mu$ that is supported in the unique  Aubry-Mather set $K$ with rotation number $\rho$. Let us explain that: as the rotations numbers $\frac{p_n}{q_n}$ are bounded, the union of these periodic orbits is bounded and then contained in some fixed compact set $K$. As the set of points having a minimizing orbit is closed, $K$ can be chosen in such a way that it is filled by minimizing orbits. Then any convergent  subsequence of $(\mu_n)$ converges to a measure $\eta$, the support of which is filled by minimizing orbits. As the rotation number is continuous on the set of points that have a minimizing orbit, the support of $\eta$ is contained in $\Mc_\rho$. As there is only one invariant Borel probability with support in $\Mc_\rho$, that is $\mu$, we conclude that $\eta=\mu$ and $(\mu_n)$ converge to $\mu$.\\

Let us now fix $\displaystyle{\sigma\in ]1, \lim_{n\rightarrow +\infty}\Rc_n[}$.  There exists $N\geq 1$ such that, for every $n\geq N$, we have $|2-{\rm Tr}(Df^{q_n}(x_n)|\geq 4\cdot\sigma^{q_n}$.\\
Because the orbit of $x_n$ is minimizing, $Df^{q_n}(x_n)$ has two real eigenvalues $\lambda_n$, $\frac{1}{\lambda_n}$ such that $\lambda_n\geq 1\geq \frac{1}{\lambda_n}$ (see the remark at the end of subsection \ref{ssAM}). We have
$$\forall n\geq N , \lambda_n+\frac{1}{\lambda_n}\geq 4.\sigma^{q_n}+2.$$
This implies that $$\liminf_{n\rightarrow +\infty}\frac{1}{q_n}\log\lambda_n\geq \log\sigma>0.$$
Hence for any $n\geq N$, the positive Lyapunov exponent of $\mu_n$ is larger that $\log\sigma$. Because the upper Lyapunov exponent depends in a upper semi-continuous way on the invariant measure (see for example \cite{Mac1}), this implies that the positive Lyapunov exponent of $\mu$ is larger than $\log\sigma$. Hence $\mu$ is hyperbolic. By Corollary \ref{corWH},  ${\rm supp}\mu$ is not contained in an essential invariant curve.

Moreover, by the results contained in \cite{Arna4}, $K$ is $C^1$-irregular $\mu$-almost everywhere.

\end{document}